\documentclass[12pt]{amsart}
\usepackage{amscd, amsmath, colonequals}
\usepackage[usenames]{xcolor}

\usepackage[hyperfootnotes=false]{hyperref}
\usepackage[usenames]{xcolor}

\def\phim{(\Phi_m^{-1})^{*}}
\def\phis{(\Phi^{-1})^{*}}

\def\blur#1{\relax}

\def\aad{AAd_{\underline U}}

\newtheorem{corollary}{Corollary}

\title[]{Convergence to normal forms of integrable PDEs}
\author{Dario Bambusi}

\address{Dipartimento di Matematica ``Federigo Enriques'',
  Universit\`a degli Studi di Milano, Via
  Saldini 50, 20133, Milano, Italia}
 \email{dario.bambusi@unimi.it}

\author{Laurent Stolovitch}
\address{Universit\'e C\^ote d'Azur, CNRS, LJAD, Parc Valrose
06108 Nice Cedex 02, France.}
\email{stolo@unice.fr}
\thanks{Research of L. Stolovitch was supported by ANR-FWF grant "ANR-14-CE34-0002-01" for the project ``Dynamics and CR geometry'' and by ANR grant ``ANR-15-CE40-0001-03'' for the project ``BEKAM''}


 \newtheorem{thm}{Theorem}[section]
 \newtheorem{theorem}{Theorem}[section]
 
 \newtheorem{cor}[thm]{Corollary}
 \newtheorem{prop}[thm]{Proposition}
 \newtheorem{lemma}[thm]{Lemma}
 \newtheorem{remark}[thm]{Remark}
 \newtheorem{definition}[thm]{Definition}

 \newcommand{\rv}{{\bf REMOVE}}

 \theoremstyle{definition}

 \renewcommand{\th}[1]{\begin{thm}\label{#1}}
 \newcommand{\eth}{\end{thm}}
 \newcommand{\co}[1]{\begin{cor}\label{#1}}
 \newcommand{\eco}{\end{cor}}
 \renewcommand{\le}[1]{\begin{lemma}\label{#1}}
 \newcommand{\ele}{\end{lemma}}
 \newcommand{\pr}[1]{\begin{prop}\label{#1}}
 \newcommand{\epr}{\end{prop}}

 \newcommand{\ga}{\begin{gather}}
 \newcommand{\ega}{\end{gather}}
 \newcommand{\gan}{\begin{gather*}}
 \newcommand{\egan}{\end{gather*}}
 \newcommand{\al}{\begin{align}}
 \newcommand{\eal}{\end{align}}
 \newcommand{\aln}{\begin{align*}}
 \newcommand{\ealn}{\end{align*}}
 \newcommand{\eq}[1]{\begin{equation}\label{#1}}
 \newcommand{\eeq}{\end{equation}}

 \newcommand{\cB}{\mathcal B}

 \newcommand{\id}{\operatorname{I}}

 \newcommand{\ve}{\vec{e}}

 \def\im{{\rm i}}
 \def\bF{{\bf F}}\def\bS{{\bf S}}
 \def\bT{{\bf T}}
 \def\bG{{\bf G}}
 
 \def\cU{{\mathcal U}}
 
 \def\cT{{\mathcal T}}
 \def\cH{{\mathcal H}}
 \def\cK{{\mathcal K}}
 \def\cX{{\mathcal X}}
 \def\bE{{\bf E}}
 \def\bR{{\bf R}}
 
 \def\bN{{\bf N}}\def\bB{{\bf B}}
 \def\sgn{{\rm sgn}}
 \def\cN{{\mathcal N}}
 \def\cNF{{\mathcal{N\kern -3pt F}}}
 \def\Z{{{\mathbb Z}^*}}
 
 \def\uno{{\rm id}}

 \newcommand{\re}[1]{(\ref{#1})}
 
 \newcommand{\rl}[1]{Lemma~\ref{#1}}

 \def\N{{\bf N}}
 \def\norma#1{\left\| #1  \right\|}
 \def\normap#1{\left\| #1  \right\|_+}

 \newcounter{pp}
 \newcommand{\bpp}{\begin{list}{$\hspace{-1em}\alph{pp})$}{\usecounter{pp}}}
 \newcommand{\epp}{\end{list}}

 \newcounter{ppp}
 \newcommand{\bppp}{\begin{list}{$\hspace{-1em}(\roman{ppp})$}{\usecounter{ppp}}}
 \newcommand{\eppp}{\end{list}}

 \def\beq{\begin{equation}}
 \def\eeq{\end{equation}}

\begin{document}
 \begin{abstract} 
 In an infinite dimensional Hilbert space we consider a family of
 commuting analytic vector fields vanishing at the origin and which
 are nonlinear perturbations of some fundamental linear vector
 fields. We prove that one can construct by the method of Poincar\'e
 normal form a local analytic coordinate transformation near the
 origin transforming the family into a normal form. The result applies
 to the KdV and NLS equations and to the Toda lattice with periodic
 boundary conditions. One gets existence of Birkhoff coordinates in a
 neighbourhood of the origin. The proof is obtained by directly
 estimating, in an iterative way, the terms of the Poincar\'e normal
 form and of the transformation to it, through a rapid convergence
 algorithm.
 \end{abstract}

 \date{\today}
  \maketitle

 \setcounter{section}{0}
 \setcounter{thm}{0}\setcounter{equation}{0}

\section{Introduction}\label{intro}

In a Hilbert space $H$, consider a family $\{X^i\}$ of (germs of) analytic vector fields defined in a
neighborhood of a common singular point, say the origin. We assume
that they are pairwise commuting with respect to the Lie
bracket. Consider the Taylor expansion $X^i=E^i+R^i$ of the fields at
the origin, with $E^i$ the linear part. It is known since
Poincar\'e, that each one of these vector fields can be transformed, by a
{\bf formal} change of variables $\hat\cT_i$ into a Poincar\'e normal
form $\hat X^i=(\hat\cT_i)_*X^i:=D\hat\cT_i(\hat\cT_i^{-1})X^i(\hat\cT_i^{-1})$. By
definition, it means that the Lie bracket $[E^i,\hat X^i]=0$
vanishes. We then say that $\hat \cT_i$ {\it normalizes} $X_i$. Since
the family is abelian, i.e. $[X^i,X^j]=0$ for al $i,j$, then one can
show that there is a single $\hat\cT$ that normalizes {\it
  simultaneously} the $X^i$'s in the sense that $[E^i, \hat\cT_*X^j]=0$, for all $i,j$.

In the same spirit, if $H$ is a symplectic {space}, one can study a family
{$\{\cH^i=\cH^i_2+ h.o.t\}$} of (germs of) analytic Hamiltonian functions which are {higher order perturbations of quadratic Hamiltonians $\cH^i_2$} and which are
pairwise commuting with respect to the Poisson bracket associated to a
symplectic form $\omega$.  The normal forms of the Hamiltonians {$\hat\cT^*\cH^i:=\cH^i\circ \hat\cT$} are then called {\it Birkhoff normal form}. We have {$\{\hat\cT^*\cH^i, \cH^j_2\}=0$} for all $i,j$ and $\hat \cT$ is a formal {\it symplectomorphism},
i.e. $\hat\cT_*\omega=\omega$.

A classical and fundamental problem in dynamics is to know under which
assumption the normalizing transformation is not only formal, but also
analytic. The motivation is to understand {\it on the normal forms}
themselves many dynamical and geometrical properties which are not
tractable directly on the original system. 
 In finite dimension, this
    problem was solved by H. R\"ussmann \cite{russmann-67} for a
    single Hamiltonian vector field and by A.D. Brjuno \cite{Bruno}
    for a single general germ of analytic vector field.  In both
    cases, one assumes that "the" formal normal form is of very
    special type, namely it has a very peculiar structure, nowaday
    called "completely integrable". For instance, in the Hamiltonian
    case, the formal Birkhoff normal form of the Hamiltonain $\cH_2+h.o.t$
    should be of the form $\hat F(\cH_2)$, with a formal power series
    $\hat F(E)=E+h.o.t$ of the single variable $E$ (a different proof
    of the same result, avoiding superconvergence has been given in
    \cite{LocatelliMele}, by developing the methods of
    \cite{LocGioZAMP}).  Still in finite dimension, J. Vey proved two
    distinct results in the same spirit. On the one hand, he
    considered in \cite{vey-ham} a family of $n$ commuting Hamiltonian
    vector fields in $\Bbb C^{2n}$, whose linear parts are linearly
    independant. On the other hand, he considered in \cite{vey-iso} a
    family of $n-1$ commuting volume preserving vector fields in $\Bbb
    C^{n}$ whose linear parts are linearly independant. In both cases,
    he proved the existence of an analytic transformation to a normal
    form of the family near the origin. In the Hamiltonian case,
    H. Ito \cite{Ito1,Ito2} impoved the results by essentially
    removing the condition of independance of the linear parts.
    N. T. Zung \cite{zung-birkhoff,zung-nf} generalized Vey's
    Hamiltonian approach by considering $m$ "linearly independent"
    vector fields having $n-m$ "functionally independant" analytic
    first integrals in $\Bbb C^n$. He proved there the convergence of
    the transformation to normal forms.  All these results have been
    unified in \cite{Stolo-ihes, Stolo-cartan} ({see also
      \cite{Stolo-asi07})} in R\"ussmann-Brjuno spirit~: it is proved
    that if the formal normal form of the family has a very peculiar
    structure (called "completely integrable"), and if the family of
    linear parts does not have "bad small divisors", then one can
    normalize analytically the family. One of the key points
    connecting the previous results with the later is that, preserving
    a structure such as a symplectic or a volume form, automatically
    implies that formal normal form of the family is "completely
    integrable".

The aim of this article is to devise such a normalizing scheme for
{\bf "complete sequences of integrable PDE's in
  involution"}. Algebraic "Hierachies of PDEs" such as defined in
\cite{magri, dickey-book} would have been the kind of objects we could
have considered but their very algebraic nature does not seem to be
suitable for our analysis.  We consider sequences of integrable
PDE's such as a family of (germs of) analytic Hamiltonian functions
$\{\cH^i\}$ in a neighborhood of a common singular point of some
suitable (infinite dimensional) Hilbert space $H$. We consider their
associated Hamiltonian vector fields $\{X^i:=X_{\cH^i}\}$ vanishing at
a common singular point, say the origin. Since the Poisson bracket
$\{\cH^i, \cH^j\}=0$ for all $i,j$, then $[X^i, X^j]=0$ for all
$i,j$. As in finite dimension, the fact that all the $X^i$'s are
symplectic implies that their formal normal form is of very special
type, namely "completely integrable" (see Definition \ref{completeint}
below). We shall show that the family of the linear parts $\{E^i\}$
at the origin, does not have "small divisors" and prove, through a rapid
convergence algorithm, that the transformation of the family to a
normal form is convergent in a neighborhood of the origin.

We also prove that our algorithm allows to construct Birkhoff
coordinates for all the vector fields $X^i$ and for all the vector
fields commuting with each one of them. We recall that Birkhoff coordinates are a
type of cartesian action angle coordinates $(x_j,y_j)$, s.t. all the
Hamiltonians of the fields $X^i$ are funcition of $x_j^2+y_j^2$ only.
We emphasize that our theorem is quite general and, as we will show,
it applies to KdV, Toda and the defocusing NLS. {Our starting point to address this problem is to consider the Lax pair \cite{lax} $\frac{dL}{dt}=[B,L]$ associated to an Integrable Pde's such as KdV. For instance, KdV equation on the circle, that is $\partial_tu-6u\partial_xu+\partial_x^3u=0$ for a function $u$ defined for $x$ on the circle $S^1$, is obtained for $L=\partial_x^2-u$, $B=-4\partial_x^3+6u\partial_x+\partial_xu$. It is known that the spectrum of $L$ is an invariant of the motion (i.e. independent of $t$) and, that the eingenvector equation turn out to be a {\it Sturm-Liouville equation} \cite{marchenko-book}. It follows that the eigenvalues can be ordered as $\lambda_0<\lambda_1\leq\lambda_2<\cdots$. As shown in \cite{garnett-trubowitz}, the sequence of {\it square of the gap lengths}, $\{(\lambda_{2n}-\lambda_{2n-1})^2\}$ forms a family of analytic first integrals commuting pairwise for a suitable Poisson bracket. Our goal is be to transform analytically and simultaneously these Hamiltonians into a Birkhoff normal form.}

We recall that a previous quite general theorem allowing to introduce
Birkhoff coordinates is due to S. Kuksin and G. Perelman
\cite{kuksin-perelman} who generalized Vey's Hamiltonian approach to
infinite dimension inspired by the scheme developed by H. Eliasson
\cite{eliasson-vey}. In the present paper we show that
Kuksin-Perelman's result can also be deduced from our Theorem
\ref{main}, in the sense that the assumptions of Kuksin-Perelman's
Theorem imply the assumptions of Theorem \ref{main}. Thus, in
particular our main result applies to all the systems for which the
assumptions of Kuksin-Perelman's Theorem hold
(\cite{kuksin-perelman,BM16,Mas18}).

We also recall that Birkhoff coordinates have been introduced
originally in PDEs by Kappeler and coworkers \cite{kapDuke,
  kappeler-poschel-book, KapHenr, KapGre, 
  kappeler-focusing-nls}. The idea of this series of
papers is to consider the square of the {spectral gaps} associated to
the Lax pair and to use them as a complete sequence of integrals of
motion in order to apply Arnold Liouville procedure\cite{Arn1} of
construction of action angle variables (which of course has to be
suitably generalized). Finally one regularizes the singularities
of such variables by introducing cartesian type coordinates, which are
the Birkhoff coordinates.

We emphasize that, althought Kuksin-Perelman's and Kappeler's approaches
are different, they are intrinsically based on the symplectic structure
and on Hamiltonian techniques.

In the present paper, we manage to {\it directly normalize} simultaneously
the family of the first integrals by a {\bf Newton scheme} (i.e. rapid
convergence scheme such as for Nash-Moser theorem
\cite{berti-corsi-procesi-NM,berti-bolle-proccesi-NM}). Furthermore,
we emphasize that our scheme is finally unrelated to the symplectic
geometry. As in finite dimension, symplectic geometry ensures that the
formal Birkhoff normal forms of all the integrals are ``completely
integrable'', i.e. are of a very special form. Such a special
form is crucial in order to estimate to solution of {\it
  nonlinear cohomological equations}.

In order to apply the algorithm in the
present infinite dimensional context, we have to face several
difficulties: the first one is to find a suitable norm to measure the
size of a family of analytic vector fields, and the second one is the
Lemma \ref{estim-nonlin-cohomo} which allow us to estimate the
``nonlinear cohomological equation'' without any small divisor
problem. The last difficulty are located in Lemmas \ref{flow} and
\ref{flow.E}, which allow us to estimate the remainder and flows
under the complete integrability assumption.

We also expect that our technique can be generalized to
the case of systems preserving other structures, e.g. a volume
form. Here we did not develop this because we are not aware of
meaningfull examples to which such a theory would apply.

We recall that it is known how to put a system in normal form up
  to some reminder in a neighbourhood of a nonresonant fixed point
  (see e.g. \cite{bam03,BG06,BDGS,bammontreal}), however the
  technique we use here is completely different from the one of these
  papers, and we do not think that the ideas of those papers, applied
  to integrable PDEs could lead to the convergence result that we
  prove here. 

Finally we remark that normal form results are often a fundemantal
starting points for studying the stability of perturbed integrable
PDEs (see e.g. \cite{kappeler-poschel-book,maspero-procesi,BKM18}).

\noindent {\it Acknowledgements.} We thank Michela Procesi for many
discussions. We acknowledge the support of Universit\`a degli Studi
di Milano and of GNFM.

\section{Main results}\label{main.sect}

\subsection{Families of normaly analytic vector fields}

Having fixed two sequences of weights $w^{(2)}_j\geq w^{(1)}_j>0$, $j\geq 1$,
we define the Hilbert spaces $H=\ell^2_{w^{(1)}}$ and $H^+:=\ell^2_{w^{(2)}}$
where $\ell^2_{w^{(n)}}$ is the Hilbert spaces of the complex
sequences $z:=\left\{z_j\right\}_{j\in\Z}$, $\Z:={\mathbb Z}-\{0\}$ s.t.
\begin{equation}
\label{nomr}
\left\|  z\right\|^2_{w^{(n)}}:=\sum_{j\in\Z}w_{|j|}^{(n)}|z_j|^2<\infty\ .
\end{equation}
In the following we will denote the norms simply by
$\left\|z\right\|:=\left\| z\right\|_{w^{(1)}}$, and
$\left\|z\right\|_+:=\left\| z\right\|_{w^{(2)}}$. Furthermore, we will
denote by ${\bf e}:=\{\vec{e_j}\}_{j\in\Z}$ the vectors with components $(\vec
e_j)_k\equiv\delta_{j,k}$, the Kronecker symbol.

Let $Q\equiv (...,q_{-k},...,q_{-1},q_1,...,q_k,...)\in \N^{\Z}$ be
an integer vector with finite support, then we write
$$
z^Q:=... z_{-k}^{q_{-k}}...z_{-1}^{q_{-1}}z_{1}^{q_{1}}.... z_{k}^{q_{k}}...,\quad |Q|:=\sum_{k\geq 1}\left(|q_{-k}|+|q_{k}|\right).
$$
We shall denote $\N_k^{\Z}$ the set of $Q\in \N^{\Z}$ with $|Q|\geq k$.

A formal vector field $X$ is a formal sum of the form
\begin{equation}
  \label{poli.for}
X(z)=\sum_{r\geq0}\sum_{i, |Q|=r} X_{Q,i}z^Q\vec{e_i}\ ,
\end{equation}
or simply 
\begin{equation}
  \label{poli.for1}
X(z)=\sum_{i, Q} X_{Q,i}z^Q\vec{e_i}\ .
\end{equation}
Two formal vector fields will be said to be equal if the corresponding
coefficients $X_{Q,i}$ coincide. 
 A formal vector field $X$ is formally conjugate
to a formal vector field $Y$, if there exists a formal vector
vector field $U$ such that
$$
Y:=(\exp U)^*X:=\sum_{k\geq 0}\frac{ad_U^k(X)}{k!},\quad ad_U(X):=[U,X],
$$
and $[.,.]$ is the commutator of vector fields. 

{ 
	\begin{definition}Let $\cU\subset H$ be a neighborhood of the origin in $H$. A formal vector field $X$ as defined by \re{poli.for1} is said to be analytic from $\cU$ to $H^+$ if the series \re{poli.for1} converges in $H^+$ uniformly for $z$ in $\cU$. The space of such vector fields will be denoted by $\cX^{\omega}(\cU, H^+)$. The space of germs at the origin of analytic vector fields with value in $H^+$ will be denoted by $\cX^{\omega}_0(H,H^+)$.
\end{definition}}

Let {$X\in \cX^{\omega}_0(H,H^+)$ be a (germ of) analytic vector field at the origin of $H$ into $H^+$} and
consider the vector field
\begin{equation}
\label{nor.an}
\underline{X}(z):=\sum_{Q,i} |X_{Q,i}|z^Q\vec{e_i}\ ,
\end{equation}
which in general is defined only on a dense subset of an open ball
$H$.
\begin{definition}
Let {$X\in \cX^{\omega}_0(H,H^+)$} be an analytic vector field
vanishing at the origin. We shall say that $X$ is normally
analytic in a ball of radius $r$ if $\underline X$ is analytic in a
ball of radius $r$ {(in $H$ with values in $H^+$)}. In this case we will write $X\in\cN_r$. We will
write $X\in\cN$ in all the cases where the value of $r$ is not important.
\end{definition}

\begin{remark}
  \label{matrix}
The above definition immediately extends to the case of applications
from $H$ to a general Banach space. In particular we will use it in 
Subsection \ref{KPsec} 
for the case where the target space is the space $B(H,H^+)$ of bounded
linear operators from $H$ to $H^+$. 
\end{remark}

{In what follows, all {\it analytic} vector fields will be considered as defined in a neighborhood (precised or not) of the origin of $H$ with values in $H^+$.}

A norm on {$\cN_r$} is given by
\begin{equation}
\label{norm.noran}
\norma{\underline X}_r:=\sup_{\norma z\leq r}\norma{\underline
  X(z)}_{+}\ .
\end{equation}

Let $X, Y$ be normaly analytic vector fields. We shall say that $Y$
dominates $X$ and we shall write $X\prec Y$, if $|X_{Q,i}|\leq
|Y_{Q,i}|$ for all indices.

\begin{remark}\label{prec}
In particular, if $X\prec Y$, then $\norma{\underline X}_r\leq
\norma{\underline Y}_r$ for any positive $r$.
\end{remark}

\begin{definition}
\label{def.fam}
A family ${\bf F}=\{F^{i}\}_{i\geq1}$ of normally analytic vector
fields will be said to be {\bf summable} if the vector field
\begin{equation}
\label{bF}
\underline{\bF }:=\sum_{i}\underline{F}^{i} \ 
\end{equation}
is normally analytic in a ball of radius $r$. In this case we will say that $\bF\in\cNF_r$.
\end{definition}

\begin{remark}
\label{multil}
Writing
$$
F^i(z)=\sum_{Q,j}F^{i}_{Q,j}
z^Q\vec{e}_j\ ,
$$
one has
\begin{equation}
\label{mul1}
\underline{\bF}(z)=\sum_{Q,j}\left(\sum_{i}\left|F^{i}_{Q,j}
\right|\right)z^Q\vec{e}_j\ ,
\end{equation}
so that, for any $r>0$, $\left\|\underline{\bF}\right\|_r$ bounds the norm of
each one of the vector fields of the family, that is 
$\left\|\underline{F^i}\right\|_r\leq \left\|\underline{\bF}\right\|_r$.
\end{remark}

\subsection{Normal forms} 
Consider the family $\bE=\left\{E^i\right\}_{i\geq 1}$ of
linear vector fields
\begin{equation}
\label{Ei}
E^i(z):=z_i\vec e_{i}-z_{-i}\vec e_{-i}\ .
\end{equation}
We will often use the notation 
\begin{equation}
\label{lambda}
E^i=\sum_{j\in \Bbb Z^*}\mu^{i}_jz_j\vec e_j\ ,
\quad \mu^i_j:=\delta^i_j-\delta^i_{-j} \ ,
\end{equation}
which is ready for the generalization to the non Hamiltonian case.

\begin{remark}
  \label{notE}
If the sequence $\frac{w_j^{(2)}}{w_j^{(1)}}\to\infty$, then the
  family $\bE$ is not summably normally analytic according to our
  definition. Indeed, the vector field $\underline\bE$ is the
  identity, which is not analytic as a map from $H$ to $H^+$.
\end{remark}

Let $\cN^{res}$ be the centralizer of the family ${\bf
  E}$, that is
$$
\cN^{res}:=\{F\in \cN \;|\; [E^i, F]=0,\quad \forall i\}
$$ By the definition of $\bE$, we have
$$
[E^i, z^Q\ve_j]=\left(\sum_{l\geq 1}q_l\mu_{l}^i - \mu_{j}^i\right)z^Q\ve_j =: \left((Q,\mu^i)-\mu_{j}^i\right)z^Q\ve_j.
$$ 
Hence, any function $F\in{\cN}^{res}$ is obtained as the
(possible infinite) linear combination of the monomials $z^Q\ve_j$ for
which $\left((Q,\mu^i)-\mu_{j}^i\right)=0$ for all $i\geq 1$.  

Let $\cN^{nres}$ be the subspace of $\cN$
generated by monomials $z^Q\ve_j$ for which
$\left((Q,\mu^i)-\mu_{j}^i\right)\neq 0$ for some $i\geq 1$.

So, any vector field $F\in\cN$ can be uniquely decomposed as
$$
F=F^{res}+F^{nres}\ ,\quad
F^{res}\in\cN^{res}\ ,\ F^{nres}\in\cN^{nres}\ .
$$ A vector field $F\in\cN^{res}$ will be called {\it resonant}, while
a vector field $F\in\cN^{nres}$ will be called {\it non
  resonant}. When speaking of the vector field $U$ which generates a
coordinate transformation we
shall say that it is {\it normalized} if $U^{res}=0$.

The same notation and terminology will be used also for families of
vector fields, and in such a case we will write $\cNF^{nres}$ for a
nonresonant family, namely a family composed by nonresonant vector fields
and similarly for $\cNF^{res}$.

\subsection{Cohomological equation} {\label{cohom}}
Let us consider the map $d_0$ which maps a homogeneous polynomial
vector field $U$ of degree $d$ to the following family of homogeneous
polynomial vector fields of degree $d$:
$$
d_0(U):=\left([E^i,U]\right)_{i\geq 1}.
$$
This map is called the cohomological operator.

A family ${\bf F}=\{F^{i}_d\}_{i\geq 1}$ of homogeneous formal
polynomial vector field of degree $d$ is called a {\bf cocycle}
with respect to the family ${\bf E}=\{E^i\}_{i\geq 1}$,if it satisfies:
\begin{equation}\label{compat}
[E^i, F^{j}_d]=[E^j, F^{i}_d],\quad i,j\geq 1.
\end{equation}
Let us write
$$
F^{j}_d=\sum_{|Q|=d,i}F^{j}_{Q,i}z^Q\ve_i.
$$ 
Therefore, equation \re{compat} reads, for all $Q\in \Bbb N^{\Bbb
  Z^*}$ and $i,j\geq 1$~: \eq{compat-monom}
\left((Q,\mu^i)-\mu_{l}^i\right)F^{j}_{Q,l}=\left((Q,\mu^j)-\mu_{l}^j\right)F^{i}_{Q,l}.
\eeq As already pointed out, any 
cocycle ${\bf F}$ can be uniquely decomposed into a sum
${\bf F}= {\bf F}^{res}+ {\bf F}^{nres}$.

\begin{lemma}\label{semisimple}
Let ${\bf F}=\{F^{i}\}_{i\geq 1}$ a formal homogeneous polynomial
vector field of degree $d$ be a non resonant cocycle (i.e satisfying
\eqref{compat}). Then, it is a coboundary, that is there exists a formal
homogeneous polynomial vector field $U$ of degree $d$ solution of the
{\bf cohomological equation} \beq \label{lin-cohomo-equ} d_0(U)={\bf
  F}, \eeq that is $[E^i,U]=F^i$, for all $i$. Furthermore, there
exists a unique normalized $U$ s.t. \eqref{lin-cohomo-equ} holds.
\end{lemma}
\begin{proof}
For each multiindex $Q\in \Bbb N^{\Bbb Z^*}_2$ and index $l\geq 1$ {such that $(F^j_{Q,l})_j\neq 0$}, 
there exists $i(Q,l)$ such that $(Q,\mu^i)-\mu_{l}^i\neq
0$. Then set $U_{Q,l}:=
\frac{F^{i(Q,l)}_{Q,l}}{(Q,\mu^{i(Q,l)})-\mu_{l}^{i(Q,l)}}$ and
$U=\sum_{Q,i}U_{Q,i}z^Q\ve_i$. Then according to \eqref{compat-monom},
we have
\begin{align*}
[E^j,U]= \sum_{Q,l}((Q,\mu^j)-\mu_{l}^j)U_{Q,l}z^Q\ve_l=
\sum_{Q,l}((Q,\mu^j)-\mu_{l}^j)
\frac{F^{i(Q,l)}_{Q,l}}{(Q,\mu^{i(Q,l)})-\mu_{l}^{i(Q,l)}}z^Q\ve_l
\\
=\sum_{Q,l}F_{Q,l}^jz^Q\ve_l.
\end{align*}

\end{proof}
\begin{definition}
  The family $\{E^i\}$ of linear vector field is said be {\bf small
    divisors free} if there exists a positive constant $c$, such that for each $Q\in \Bbb N^{\Bbb Z^*}_2$ and $j\geq 1$, there is $i(Q,j)\geq 1$ such that $|(Q,\mu^{i(Q,j)})-\mu^{i(Q,j)}_j|>c^{-1}$.
\end{definition}
\begin{remark}\label{rem-sd}
If the  family ${\bf E}$ is small divisor free, then for any summable normally analytic cocycle ${\bf F}$, the unique normalized solution $U$ to \eqref{lin-cohomo-equ} is normally analytic and satisfies, for some $r>0$
$$
\norma{\underline U}_r\leq c \norma{\underline {\bf F}}_r
$$
for some constant $c$.
\end{remark}
\begin{definition}
  \label{nf.1}
A formal vector field $X$ (resp. a family ${\bf X}=\{X^{j}\}_{i\geq
  1}$) is said be a {\bf normal form} with respect to ${\bf E}$ if
$[E^i,{X}]=0$ for all $i\geq 1$ (resp. $[E^i,X^{j}]=0$ for all
$i,j\geq 1$).
\end{definition}

\begin{remark}
  \label{rem.norfor}
The above definition is taylor made for the case of Hamiltonian vector
fields. In the case of vector fields preserving different structures
the definition has to be modified following \cite{Stolo-cartan}.
\end{remark}

\begin{definition} \label{completeint}
An analytic (resp. formal) normal form $X$ is said to be {\bf completely integrable}
if it can be written as $X=\sum_{j\geq 1}a_j E^j$ where $a_j$ are
normally analytic (resp. formal) functions, invariants w.r.t ${\bf E}=\{E^i\}_{i\geq
  1}$, i.e. $E^i(a_j)=0$, for all $i,j\geq 1$.
\end{definition}
\begin{definition}
A family of formal vector fields is said to be {\bf formally completely
integrable} if it is formally conjugate 
to a completely integrable formal normal form.
\end{definition}

\begin{lemma}
 A formal transformation of the form $\exp(U)$ with $U=\sum_{k\geq 2}U_k$, where $U_k$ is a homogeneous formal polynomial of degree $k$ commuting with each $E_i$, $i\geq 1$, conjugates a formal normal
    form of a family of formal vector fields to another normal form.

If one of the formal normal forms 
is completely integrable, so are all the other
normal forms.
\end{lemma}

\begin{proof}
First of all, if $\{Z^j\}$ and $\{(\exp U)_*Z^j\}$ are normal forms then $[U,E^i]=0$ for all $i$. Indeed, we have
$$
(\exp U)_*Z^j=Z^j+[U,Z^j]+\frac{1}{2}[U, [U,Z^j]]+\cdots
$$
Taking the bracket with $E^i$ and using Jacobi identity, we obtain
$$
[E^i,(\exp U)_*Z^j]=-[Z^j,[E^i,U]]+\frac{1}{2}[E^i,[U, [U,Z^j]]]+\cdots
$$
Let $U_{d_0}$ be lowest order term in the Taylor expansion of $U$ at origin. Then, one has
$[E^j,[E^i,U_{d_0}]]=0$ for all $i,j$. Hence, $[E^i,U_{d_0}]$ belongs to both the range and the kernel of the semi-simple map $[E_i,.]$. Hence, $[E^i,U_{d_0}]=0$ for all $i$. 
Hence, the lowest order term of $[Z^j,[E^i,U]]$ is $[E^j,[E^i,U_{d_0+1}]]$. On the other hand, since the bracket of resonant vector fields is still resonant, we have, for $k\geq 2$, 
$[E^i,ad_U^k(Z^j)]=[E^i,ad_{U}^{k-1}([U-U_{d_0},Z^j])$ which is of order $\geq d_0+1+(k-1)(d_0-1)\geq 2d_0>d_0+1$. Hence, the lowest order term of $[E^i,(\exp U)_*Z^j]$ is $[E^j,[E^i,U_{d_0+1}]]=0$ and we proceed by induction on the order.

Assume that the the family $\{Z^j\}$ is completely integrable. Transform it to another normal form ${\tilde Z^j}$ by a transformation $\exp U$. According to the first point, $U$ commutes with each $E^i$. Hence, it commutes with each $Z^i$ since $[U,\sum_ja_{i,j}E^j]=\sum_ja_{i,j}[U,E^j]+U(a_{i,j})E^j$. On the other hand, $E^k(U(a_{i,j}))= [E^k,U](a_{i,j})=0$ for all $k$. So that 
$$
(\exp U)_*Z^i=\sum_j(\sum_k \frac{U^k(a_{i,j})}{k!}) E^j.
$$
\end{proof}
The main result of this paper is the following theorem.

\begin{theorem}
\label{main}
Consider a family of analytic vector fields of the form
\begin{equation}
\label{main1}
X^{i}=E^i+F^i\ ,\quad i\geq1\ .
\end{equation}
Assume that
\begin{itemize}
\item[0.] the family of linear vector fields $\bE$ is small divisor free.
\item[1.] $\bF\equiv\{F^i \}_{i\geq 1}$ is a summable family of
  normally analytic vector fields
\item[2.] there exists $c_0$ s.t. for $r_0$ small enough, one has
  $\norma{\underline{\bF}}_{r_0}\leq c_0r_0^2$
\item[3.] $[X^{i};X^j]\equiv 0$, $\forall i,j$.
\item[4.] ${\bf
  X}$ is formally completely integrable.
\end{itemize}
Then there exist constants $r_*>0$, $c_2,c_3$, a neighborhood
$\cU\supset B_{r_*}$ of the origin and an analytic coordinate
transformation $\cT:\cU\to H$ s.t.
\begin{equation}
\label{maine}
\cT_*X^i=E^i+N^i\ ,\quad \forall i\geq 1 \ ,
\end{equation} 
where $\cNF_{r_*}\ni \bN \equiv\left\{ N^i\right\}_{i\geq 1}$ is a completely
integrable normal form. 

{Furthermore}, $\forall r<r_*$ the following estimates hold:
\begin{itemize}
\item[i.] $\norma{\bN}_{r}\leq c_2r^2$,
\item[ii.]  $\sup_{\norma{z}\leq r}\normap{z-\cT(z)}\leq c_3r^2$.
\end{itemize}
\end{theorem}
\begin{remark}
\label{stru}
From the proof it is clear that if one endows the Hilbert space by the
symplectic structure $\im dz_{-k}\wedge dz_{k}$ and the vector fields
$F^i$ are Hamiltonian for any $i$, then the transformation $\cT$ is
canonical. Here we did not assume the fields $X^i$ to be
  Hamiltonian. In the Hamiltonian case Assumption 4 would be
  automatic. 

We expect that the result can be extended also to other preserved
structure, like volume in phase space, but the present proof rely on
the structure of the family $\bE$.  
\end{remark}

\begin{remark}
The Hilbert spaces considered can chosen to be more general. For instance, it could be be spaces of sequences indexed over $\Bbb Z^d\setminus\{0\}$.
\end{remark}

We are now going to give a more precise statement for the Hamiltonian
case, showing in particular that the transformation $\cT$ introduces
Birkhoff coordinates for the integrable Hierarchy associated to the
fields $\{X^i\}$.

Thus, in the space $H$, we introduce the symplectic form $\im dz_{-k}\wedge
dz_k$. Given an analytic function $\cH\in C^\omega(H, \Bbb R)$, we define the
corresponding Hamiltonian vector field $X_{\cH}$ as the vector field with
$k$-th component
$$
(X_{\cH})_k(z):=-i(\sgn k) \frac{\partial \cH}{\partial z_{-k}}\ .
$$ Given also a second function $\cK \in C^\omega(H, \Bbb R)$ we define
their Poisson Bracket by
$$
\{\cH,\cK \}(z):= d\cH(z) X_{\cK}(z)\ .
$$ It is well known that such a quantity can fail to be well defined,
nevertheless in all the cases we will consider it will be well
defined.

Consider now a sequence of analytic Hamiltonians $\cH^i$ of the form
$$
\cH^i=\cH^i_2+\cK^i\ ,\quad \cH^i_2:=z_iz_{-i}\ ,
$$
and $\cK^i$ having a zero of order at least $3$ at the origin.

\begin{corollary}
  \label{main.1}
Assume that the vector fields $X^i:=X_{\cH^i}$ fulfill the assumtpions
of Theorem \ref{main}, then the coordinate transformation $\cT$ is
canonical. Furthermore, given any analytic Hamiltonian $\cH$ with a
zero of order 2 at the orgin, such that
\begin{equation}
  \label{commu.h}
\{\cH,\cH^i\}\equiv 0\ ,\quad \forall i\geq 1\ ,
\end{equation}
one has that $\cH\circ\cT^{{-1}}$ is a function of
  $\{(z_jz_{-j})\}_{j\geq 1}$ only.
\end{corollary}
\proof The proof follows the proof of Corollary 2.13 of
\cite{BM16}. First, it is clear that $E^j$ is the Hamiltonian vector
field of $\cH^j_2$. Denote $\tilde \cH:=\cH\circ\cT^{{-1}}$, thus, from the
property that $\cT_*X^j$ is in normal form one has that
\begin{equation}
  \label{comm.rel.2}
[\cT_*X_{\cH},E^j]=X_{\{\tilde \cH,\cH^j_2  \}}=0\ ,
\end{equation}
from which $\{\tilde \cH,\cH^j_2 \}=c^j$. However, since both
$\tilde \cH$ and $\cH_2$ have a zero of order 2 at the orgin, the
constants must vanish. Expand now $\tilde\cH$ in Taylor
series, one has
$$ \tilde\cH (z) = \sum_{\substack{r \geq 2,
    \\ |\alpha| + |\beta| =r }} H^r_{\alpha, \beta} z_+^\alpha 
z_-^\beta\ ,
$$
where we denoted $z_+:=\{z_j\}_{j\geq 1}$ and $z_-:=\{z_{-j}\}_{j\geq
  1}$. Then equation \eqref{comm.rel.2} implies that in each term of the
summation $\alpha = \beta$, therefore $\tilde \cH$ is a
function of $z_jz_{-j}$ only.  \qed

\subsection{Kuksin-Perelman's Theorem}\label{KPsec}

In this section we recall the Vey type theorem obtained by Kuksin
and Perelman in \cite{kuksin-perelman} (see also \cite{BM16,Mas18}) and prove that it
can be obtained as a corollary of Theorem \ref{main}.

We come to the assumptions of the Kuksin-Perelman's Theorem.

Consider an analytic map $\Psi$ of the form
\begin{equation}
\label{KP0}
\Psi=\uno+G\ ,
\end{equation} 
with $G\in\cN_R$ (with some $R>0$) having a zero of second order at
the origin. For $j>0$ consider also the functions
$I_j(z):=\Psi_j(z)\Psi_{-j}(z)$ and the Hamiltonian vector fields
$X^j:=X_{I_j}$.

Assume that the following Hypotheses hold:
\begin{itemize}
\item[(KP1)] The functions $I_j(z)$ pairwise
  commute, namely $\left\{I_j;I_k\right\}\equiv 0$ forall $j,k\geq 1$. 
\item[(KP2)] the maps $\underline{dG} $ and $\underline{dG^*}$ are
  analytic as maps from $B_R$ to $B(H,H^+)$.
\end{itemize}

\begin{theorem}[Kuksin-Perelman]
  \label{KP1}
  Assume that (KP1) and (KP2) hold, then the same conclusions of
  Theorem \ref{main} and Corollary \ref{main.1} hold.
\end{theorem}
\proof It is enough to show that the assumptions (KP1-KP2) imply the
assumptions of Theorem \ref{main} with the fields $X^j:=X_{I_j}
$. First remark that assumption 3 of Theorem \ref{main} follows from
(KP1), while assumption 4 follows from the fact that the fields
$X^i$ are Hamiltonian. Assumption 0 follows from the structure
\eqref{KP0} of the function $\Psi$. 

In order to verify assumptions 1, compute explicitely the
components of the vector fields $X_{I_l}$. For $k\geq 1$ its $k-th$
component is given by
\begin{align}
\label{2.18}
\im \left(X_{I_l}\right)_k =z_k\delta_{lk}+G^k\delta_{lk}+z_l\frac{\partial
  G^{-l}}{\partial z_{-k}}+G^l\frac{\partial
  G^{-l}}{\partial z_{-k}}+z_{-l}\frac{\partial
  G^{l}}{\partial z_{-k}}+G^{-l}\frac{\partial
  G^{l}}{\partial z_{-k}}\ ;
\end{align}
the first term contribute to $E^l$, while all the other ones contribute to
$F^l$.  From \eqref{2.18} we have that the $k$-th component of
$\underline{\bF}$, ($k\geq 1$) is given by
\begin{align}
\underline{G^k}+\sum_{l\geq 1} z_l\frac{\partial
  \underline{G^{-l}}}{\partial z_{-k}}
+\sum_{l\geq 1}\underline{G^l}\frac{\partial
 \underline{ G^{-l}}}{\partial z_{-k}}
+\sum_{l\geq 1}z_{-l}\frac{\partial
  \underline{G^{l}}}{\partial z_{-k}}
+\sum_{l\geq 1}\underline{G^{-l}}
\frac{\partial
  \underline{G^{l}}}{\partial z_{-k}}\ .
\end{align}
We have to show that each one of the terms of this expression define
the $k$-th component of an analytic vector field. For the first term
this is a trivial consequence of the fact that $G\in\cN$. Consider the
second term. In order to see that it is analytic we write it in terms
of $\underline{dG^*}$. To this end define the involution
$(Iz)_k:=z_{-k}$ and the truncation operator $(Tz)_k=z_k$ if $k\geq
1$ and zero otherwise. Then the second term of the above expression
is the $-k$-th component of $\underline{dG^*}(ITz)$, which belongs to
$\cN$ by assumption (KP1). All the other terms can be dealt with in
the same way geting that assumption 1 is fulfilled. Assumption 2 is a
direct consequence of the fact that $G$ has a zero of order 2 at the
origin.  
\qed

\section{Proof of the main theorem}\label{proof}

\subsection{Nonlinear cohomological equation}

Assume the abelian family ${\bf X}=\{X^{i}\}$ is normalized up to order $m=2^k$:
$$
X^i=E^i+N^i_{\leq m}+R^i_{\geq m+1}.
$$ where $\textbf{N}_{\leq m}\in\cNF_{R_m}$ is a completely
  integrable normal form of degree $m$; we shall write $NF^i_{\leq
    m}=E^i+N^i_{\leq m}=\sum_{j\geq 1}(\delta_{i,j}+a_{i,j})E^j$ where
  $a_{i,j}$ are polynomials of degree $\leq m-1$ that are common first
  integrals of $\bf E$.
Let us Taylor expand $R^i_{\geq m+1}=B^i_{\leq 2m}+ \tilde R^i_{\geq
  2m+1}$ up to degree $2m$. We shall (mostly) omit the dependence on $m$ in this
section. 
  Since
$X^i$ and $X^j$ are pairwise commuting, then
$$
0=[X^i,X^j]=[NF^i,B^j]-[NF^j,B_i]+[NF^i,\tilde R_{\geq 2m+1}^j]-[NF^j,\tilde R_{\geq 2m+1}^i]+[\tilde R_{\geq 2m+1}^i,\tilde R_{\geq 2m+1}^j].
$$
Therefore, the truncation at degree $\leq 2m$ gives
\begin{equation}
  \label{coc.1}
0=J^{2m}([NF^i,B^j]-[NF^j,B^i]),
\end{equation}
where $J^{2m}(V)$ denotes the $2m$-jet of $V$.
\begin{lemma}
  \label{coc}
  Let $\bB$ be a nonresonant family and $\bN$ a completely integrable
  normal form. Assume that they fulfill \eqref{coc.1}. Then there exists a unique $U$ normalized (i.e. no resonant term in expansion) such that for all $j$ one has $J^{2m}([NF^j, U])=B^j$.
\end{lemma}

\begin{proof}

We give here a direct proof although a more conceptual proof involving
spectral sequences can be found in \cite{Stolo-ihes}[proposition
  7.1.1] for the finite dimensional case. For any integer $m+1\leq k \leq 2m$, the homogenous polynomial of degree $k$ of eq.\eqref{coc.1} is
\begin{equation}\label{nonlin-compt}
\sum_{p=1}^{k-m}{[NF_p^i, B^j_{k-p+1}]}=\sum_{p=1}^{k-m}{[NF_p^j, B^i_{k-p+1}]}
\end{equation}

Let us prove, by induction on the integer $m+1\leq k \leq 2m$, that there exists a unique 
normalized polynomial $V_k$ homogeneous of degree $k$, such that 
\begin{equation}\label{result-rec}
\forall 1\leq i,\;[E^i,V_k]=B_k^i+\sum_{p=2}^{k-m}{[V_{k-p+1},NF_p^i]},
\end{equation}
that is $J^k\left([NF^i,\sum_{p=m+1}^kV_p]\right)=J^k(B^i)$.

For $k=m+1$, the equation $(\ref{nonlin-compt})$ leads to 
$[E^i,Z^j_{m+1}]=[E^j,Z^i_{m+1}]$.
According to the Lemma $\ref{semisimple}$, there exists a unique normalized
$V_{m+1}$ homogeneous of degree $m+1$ such that, for all $1\leq i$,
$[E^i,V_{m+1}]=Z^i_{m+1}$.

Let us assume that the result holds for all integers $q<k$. Let $2\leq
p \leq k-m$ be an integer, then $m+1\leq k-p+1 < k$. Let us first
recall that, by assumptions, $[NF_k^i,NF_{k'}^j]=0$ for all integers
$1\leq i,j$ and $1\leq k,k'$.

Thus, by Jacobi Identity, we have
\begin{eqnarray*}
[NF_p^i,[E^j,V_{k-p+1}] ] & = & -[E^j,[V_{k-p+1},NF_p^i] ]\\\;
[NF_p^i,[V_{k-p-q+2},NF_q^{j}]] & = & -[NF_q^j,[NF_p^i,V_{k-p-q+2}] ]\quad\forall\; 2\leq q\leq k-p+1-m\\
\end{eqnarray*}
With these remarks as well as $(\ref{result-rec})$, it follows, by induction, that 
\begin{eqnarray*}
[NF_p^i, Z^j_{k-p+1}] & = & \left[NF_p^i,[E^j,V_{k-p+1}]-\sum_{q=2}^{k-p+1-m}{[V_{k-p-q+2},NF_q^j]}\right]\\
& = & -[E^j,[V_{k-p+1},NF_p^i]]+\sum_{q=2}^{k-p+1-m}{[NF_q^j,[NF_p^i,V_{k-p-q+2}]]}\\
\end{eqnarray*}
Since $[NF^j_q,[NF^i_p, V]]=[NF^i_p,[NF^j_q, V]]$, then exchanging $j$ and $i$ leads to 
$$
[NF_p^i, Z^j_{k-p+1}]+[E^j,[V_{k-p+1},NF_p^i]] = [NF_p^j, Z^i_{k-p+1}]+[E^i,[V_{k-p+1},NF_p^j]]
$$
Summing over $2\leq p\leq k$ and using the compatibility condition $(\ref{nonlin-compt})$ leads to
$$
\left[E^j,Z^i_k+\sum_{p=2}^{k-m}[V_{k-p+1},NF_p^i]\right]= \left[E^i,Z^j_k+\sum_{p=2}^{k-m}[V_{k-p+1},NF_p^j]\right]
$$
But, the same argument as in the proof of the first point of this proposition will show that, 
$\{\sum_{p=2}^{k-m}[V_{k-p+1},NF_p^i]\}$ is a non-resonant family of homogenoues vector fileds of degree $k$.
Therefore, according to Lemma $\ref{semisimple}$, there exists a
unique normalized $V_k$ such that, for all $i\geq 1$, 
$$
[E^i, V_k]= Z^i_k+\sum_{p=2}^{k-m}{[V_{k-p+1},NF_p^i]},
$$
which ends the proof of the induction and the proposition.

\end{proof}
Let us construct and estimate the unique nonresonant solution
$U$ (i.e. with $U^{res}\equiv 0$), of order $\geq m+1$ and degree
$\leq 2m$ of the \emph{nonlinear cohomological equation}, namely 
\beq\label{nonlin-cohom}
J^{2m}([NF^i, U])=B^i_{nres},
\eeq
where $B^i_{nres}$ denotes the nonresonant projection of $B^i$. 

\begin{lemma}\label{estim-nonlin-cohomo}
Assume that, for all $\|z\|\leq r$, for all $v\in H$
$$ \|D\underline{{\bf N}}(z).v\|_+\leq \frac{1}{2}\|v\|,\quad
\|\underline{\bf R}_{\geq m+1}\|_{r}\leq \epsilon.
$$ Then \re{nonlin-cohom} has a unique nonresonant solution $U$
which satisfies
$$
\|\underline{U}\|_{r}\leq 4\epsilon.
$$
\end{lemma}
\begin{proof}
Let us write \re{nonlin-cohom} as
\beq\label{nonlin-cohomb}
[NF^i, U]=B^i_{nres}+Z^i_{\geq 2m+1}=:F^i
\eeq
where $Z^i_{\geq 2m+1}:=J^{2m}[NF^i, U]-[NF^i, U]$.

Let $\lambda_i^{(d)}$ be an eigenvalue of the operator $[E^i,\cdot]$
in the space of formal homogeneous polynomial vector fields of
degree $d$. Let $h_{i,\lambda_i^{(d)}}$ be the associated eigenspace.

\begin{remark}
\label{ham.case}
Due to the definition of the family $\bE$, we have
$\lambda_i^{(d)}=q_i-q_{-i}-s_{i}$ where $q_i$ denotes the $i$th
component of a multiindex $Q=(\cdots q_{-i-1},
q_{-i},\cdots,q_{-1},q_{1},\cdots q_{i}, q_{i+1},\cdots)$ with modulus
$d$ and $s_i$ is 1,-1 or 0. Indeed, we have $[E^i, x^Q\ve_k]=
(q_i-q_{-i}-s_{i})x^Q\ve_k$ with $s_i=1$ if $k=i$, $s_i=-1$ if
$k=-i$ and $s_i=0$ otherwise. 
\end{remark}

 Let $\lambda^{(d)}=(\lambda_1^{(d)}, \lambda_2^{(d)},\ldots)$ be a
 collection of such eigenvalues. We shall say that $\lambda^{(d)}$ is
 a generalized eigenvalue of degree $d$. If $\lambda^{(d)}\neq 0$,
 then only a finite number of its components are non zero.  Let us denote
 $\text{Supp}(\lambda^{(d)})$, the support of $\lambda^{(d)}$, that is
 the set of indexes $j$ such that $\lambda^{(d)}_j\neq 0$.
From now on, we shall write $\lambda$ for $\lambda^{(d)}$, if there is
no
confusion. 

We remark that, given $U\in \cap_{i\geq 1}h_{i,\lambda_i^{(d)}}$, and
any function $a$ which is a common first integral of the family $E^i$,
namely s.t. $E^k(a)=0$, $\forall k$, one has 
$$
[E^i,aU]=\lambda_i^{(d)} aU\ ,
$$
thus it is convenient to denote
\begin{equation}
\label{autogen}
H_{\lambda^{(d)}}:=\left\{ U\in\cN\ :\ [E^i,U]=\lambda_i^{(d)} U
\right\} \ .
\end{equation}

We now show that $[NF^i,.]$ leaves invariant $H_{\lambda}$ (where we
omitted the index $d$ from $\lambda$). We have 
\beq
[NF^i, U]= [E^i,U]+\sum_{j\geq 1}a_{i,j}[E^j,U]+U(a_{i,j})E^j.\label{exp-nonlin-cohom}
\eeq
Here, $U(a_{i,j})$ denotes the Lie derivative of $a_{i,j}$ along $U$.
Since the $E^i$'s are pairwise commuting and since the $a_{i,j}$'s are first integrals of $\bf E$, we have
\begin{eqnarray*}
\left[E^l, [E^i,U]+\sum_{j\geq 1}a_{i,j}[E^j,U]\right]&= &[E^i, [E^l,U]+\sum_{j\geq 1}a_{i,j}[E^j,[E^l,U]]\\
&=& \lambda_l  \left([E^i,U]+\sum_{j\geq 1}a_{i,j}[E^j,U]\right).
\end{eqnarray*}
On the other hand, we have 
\begin{eqnarray*}
\left[E^l, \sum_{j\geq 1}U(a_{i,j})E^j \right] &=& \sum_{j\geq 1}E^l(U(a_{i,j}))E^j=  \sum_{j\geq 1}[E^l, U](a_{i,j})E^j\\
&=& \lambda_l\sum_{j\geq 1} U(a_{i,j})E^j
\end{eqnarray*}
From which the invariance of $H_\lambda$ follows.

Let $U_{\lambda}$ (resp. $F^i_{\lambda}$) be the projection onto $H_{\lambda}$ of $U$ (resp. $F^i$). Therefore, the projection onto $H_{\lambda}$ of equation \re{nonlin-cohomb} reads
\beq\label{proj-nonlin-cohom}
[NF^i, U_{\lambda}]=F^i_{\lambda}.
\eeq
Using \re{exp-nonlin-cohom}, this equation reads
$$
\left(\lambda_i+\sum_{j\geq 1}a_{i,j}\lambda_j\right)U_{\lambda}+ \sum_{j\geq 1}U_{\lambda}(a_{i,j})E^j= F^i_{\lambda}.
$$
Let $\epsilon_i$ be the sign of $\lambda_i$, if $i\in \text{Supp}(\lambda)$. Let us multiply the $i$th-equation by $\epsilon_i$ and then let us sum up over $i\in \text{Supp}(\lambda)$. We obtain
\begin{equation}
  \label{3.22a}
\left(|\lambda|+\sum_{j\geq 1}\sum_{i\in
  \text{Supp}(\lambda)}\epsilon_i a_{i,j}\lambda_j\right)U_{\lambda}+
\sum_{j\geq 1}\sum_{i\in
  \text{Supp}(\lambda)}\epsilon_iU_{\lambda}(a_{i,j})E^j= \sum_{i\in
  \text{Supp}(\lambda)}\epsilon_iF^i_{\lambda}=:\tilde F_{\lambda}. 
\end{equation}
Let us define
\beq\label{blamdda} b_{\lambda}:= |\lambda|+\sum_{j\in
  \text{Supp}(\lambda)}\sum_{i\in \text{Supp}(\lambda)}\epsilon_i
a_{i,j}\lambda_j=:|\lambda|+c_{\lambda}.
\eeq
Remark that it is an analytic function whose value at $0$ is
$|\lambda|$; furthermore one has $E^j(b_\lambda)=0$, $\forall j$. Let
us consider the operator $$ P_{\lambda}~: U_{\lambda}\mapsto
\sum_{j\geq 1}\sum_{i\in
  \text{Supp}(\lambda)}\epsilon_iU_{\lambda}(a_{i,j})E^j.
$$
We have $P_{\lambda}^2=0$. Indeed, since the $a_{i,j}$ are first integrals of $\bf E$, we have 
\begin{eqnarray*}
P_{\lambda}(P_{\lambda}(U_{\lambda}))& =& \sum_{j\geq 1}\sum_{i\in \text{Supp}(\lambda)}\epsilon_i P_{\lambda}(U_{\lambda})(a_{i,j})E^j\\
&=& \sum_{j\geq 1}\sum_{k\geq 1}\sum_{i\in \text{Supp}(\lambda)}\epsilon_i\sum_{l\in \text{Supp}(\lambda)}\epsilon_l U_{\lambda}(a_{l,k})E^k(a_{i,j})E^j\\
&=&0.
\end{eqnarray*}
Similarly one has $P_\lambda(P_{\lambda}(./b_\lambda))=0$.  As a
consequence, the nonresonant solution of equation
\re{3.22a} is
\beq\label{sol-proj-nonlin-cohom}
U_{\lambda}= (I-\frac{1}{b_{\lambda}}P_{\lambda})\left(\frac{\tilde
  F_{\lambda}}{b_{\lambda}}\right).
\eeq
Summing up over the set of
generalized eigenvalues $\lambda$ of degree $m+1\leq d\leq 2m$, and
applying $J^{2m}$ we
obtain
\begin{equation}
  \label{Um} U=J^{2m}\left( \sum_{\lambda}\frac{\tilde
  F_{\lambda}}{b_{\lambda}}-\sum_{\lambda}\frac{1}{b_{\lambda}}P_{\lambda}\left(\frac{\tilde
    F_{\lambda}}{b_{\lambda}}\right)
  \right).
\end{equation}
Since $U$ is of
degree $\leq 2m$, we can substitute $B_\lambda$ to $F_\lambda$, thus we are led to the final definition of $U$, namely
\begin{equation}
  \label{Um.1}
U=J^{2m}\left( \sum_{\lambda}\frac{\tilde
  B_{\lambda}}{b_{\lambda}}-\sum_{\lambda}\frac{1}{b_{\lambda}}P_{\lambda}\left(\frac{\tilde
    B_{\lambda}}{b_{\lambda}}\right)
  \right),
\end{equation}
where
 $\tilde B_{\lambda}:=\sum_{i\in
  \text{Supp}(\lambda)}\epsilon_iB^i_{\lambda}$.  We now estimate such
a quantity.
Remark first that one has
$$ \underline{\left(\frac{1}{b_\lambda}\right)}=\underline{\left(\frac{1}{|\lambda|-c_\lambda}\right)}\preceq
\frac{1}{|\lambda|}\sum_{k\geq 0}\left(\frac{\underline
  {c_\lambda}}{|\lambda|}\right)^k\preceq
\frac{1}{|\lambda|-\underline{c_\lambda}} ,
$$
so that we have
$$
\underline{\sum_{\lambda}\frac{\tilde B_{\lambda}}{b_{\lambda}}} \prec  \sum_{\lambda}\frac{\sum_{i\in \text{Supp}(\lambda)}\underline{B^i_{\lambda}}}{|\lambda|-\underline{c_{\lambda}}}
$$ On the other hand, given an orthonormal basis ${\bf e}$ of $H^+$, a sequence
$\{G_{\lambda}\}$ of vectors with nonnegative coordinates on ${\bf e}$
and a bounded sequence $\{g_i\}$ of nonnegative numbers, we have
$$ \left\|\sum_{\lambda}g_{\lambda}G_{\lambda}\right\|_+^2=
\sum_{{\lambda},{\lambda}'}g_{\lambda}g_{\lambda'}(G_{\lambda},G_{\lambda'})_+\leq
(\sup_{{\lambda},{\lambda'}}g_{\lambda}g_{\lambda'})\left\|\sum_{\lambda}G_{\lambda}\right\|_+^2\leq
(\sup_{{\lambda}}g_{\lambda})^2\left\|\sum_{\lambda}G_{\lambda}\right\|_+^2.
$$
Evaluating at a point near the origin in the
domain, we can apply this with
$g_{\lambda}=\frac{1}{|\lambda|-\underline{c_{\lambda}}} $ and
$G_{\lambda}=\sum_{i\in
  \text{Supp}(\lambda)}\underline{B^i_{\lambda}}$ 
	Hence, we obtain \beq\label{maj-Bb} \left\|
\sum_{\lambda}\frac{\sum_{i\in
    \text{Supp}(\lambda)}\underline{B^i_{\lambda}}}{|\lambda|-\underline{c_{\lambda}}}\right\|_+\leq
\sup_{\lambda}\left|\frac{1}{|\lambda|-\underline{c_{\lambda}}}\right|\left\|
\sum_{\lambda}\sum_{i\in
  \text{Supp}(\lambda)}\underline{B^i_{\lambda}}\right\|_+\leq
\sup_{\lambda}\left|\frac{1}{|\lambda|-\underline{c_{\lambda}}}\right|\|\underline{{\bf
    B}}\|_+.  \eeq In order to estimate $\underline{c_\lambda}$,
remark first that according to \re{blamdda}, we have \beq\label{maj-c}
\underline{c_{\lambda}}\prec \sum_{j\in
  \text{Supp}(\lambda)}\sum_{i\in
  \text{Supp}(\lambda)}\underline{a_{i,j}}|\lambda_j|\prec \sum_{j\in
  \text{Supp}(\lambda)}|\lambda_j|\left(\sum_i
\underline{a_{i,j}}\right).  \eeq To estimate
$\beta_j:=\sum_{i}\underline{a_{ij}}$, we proceed as
follows. According to \re{lambda}, we have
$\underline{N^i}=\sum_{j\in\Bbb Z^*}\underline{a_{i,j}}z_j\ve_j$ so
that $\underline{{\bf N}}=\sum_{j\in\Bbb
  Z^*}\left(\sum_i\underline{a_{i,j}}\right)z_j \ve_j=
\sum_{j}\beta_jz_j\ve_j$. Hence, we have
$$
\frac{\partial\underline{{\bf N}}}{\partial z_k}=\sum_{j\in\Bbb Z^*}\frac{\partial \beta_j}{\partial z_k} z_j\ve_j+\beta_k\ve_k.
$$
Since the previous equality involves only vectors with nonnegative coefficients, we have
\begin{eqnarray}
\beta_k e_k &\prec& \frac{\partial\underline{{\bf N}}}{\partial
  z_k}\label{apartialnf}
\\
\sum_{j\in\Bbb Z^*}\frac{\partial \beta_j}{\partial z_k} z_je_j&\prec & \frac{\partial\underline{{\bf N}}}{\partial z_k}\label{partialapartialnf}
\end{eqnarray}
So, $\forall v\in H$ and for all $\|z\|\leq r$, we have
\begin{align*}
\left(\frac{1}{2}\norma v \right)^2=\frac{1}{4}\sum_{k}w_k^1{(1)}\left|
v_k\right|^2 \geq \normap{D\underline{\bN}(z)v}^2=
\normap{\sum_{k}\frac{\partial\underline{\bN}}{\partial z_k}v_k }^2 \\ \geq
\normap{\sum_{k}v_k\beta_k\ve_k}^2=\sum_{k}w_k^{(2)}\beta_k^2v_k^2 =
\sum_{k}w_k^{(1)}\frac{w_k^{(2)}}{w_k^{(1)}}\beta_k^2v_k^2\ .
\end{align*}
Taking $v:=v_k\ve_k=1/\sqrt{w_k^{(1)}}\ve_k$, which has norm 1, one gets
$$
\frac{1}{4}\geq \frac{w_l^{(2)}}{w_l^{(1)}}\beta_l^2\geq \beta_l^2 =
\left(\sum_{i}\underline a_{il}\right)^2\ . 
$$
Inserting in \eqref{maj-c} one gets
$$
\left|\underline c_\lambda\right|\leq
\left|\lambda\right|\frac{1}{2}\ ,
$$
hence
$$
\left|\frac{1}{|\lambda|-\underline{c_{\lambda}(z)}}\right|\leq \frac{2}{|\lambda|}.
$$

Since the familly $\bE$ is small divisor free, then we always have $1\leq |\lambda|$ (we have set $c=1$ for simplicity), then by \re{maj-Bb}
\beq\label{estim1}
\sup_{\|z\|\leq r_m}\left\|\underline{\sum_{\lambda}\frac{\tilde B_{\lambda}}{b_{\lambda}}}\right\|\leq \sup_{\lambda}\frac{2\epsilon}{|\lambda|}\leq 2\epsilon
\eeq
as soon as $\|\underline{{\bf B}}\|_+\leq \epsilon$.
On the other hand, we have
$$
\sum_{\lambda}\frac{1}{b_{\lambda}}P_{\lambda}\left(\frac{\tilde B_{\lambda}}{b_{\lambda}}\right)\prec \sum_{\lambda}\frac{1}{(|\lambda|-\underline{c_{\lambda}})^2}\sum_{j\geq 1}\sum_{i\in \text{Supp}(\lambda)}D\underline{a_{i,j}}\left(\sum_{l\in \text{Supp}(\lambda)}\underline{B^l_{\lambda}}\right)\underline{E^j}
$$
According to \re{partialapartialnf}, we have
$$
\sum_{j\geq 1}\sum_{i\in \text{Supp}(\lambda)}D\underline{a_{i,j}}\left(\sum_{l\in \text{Supp}(\lambda)}\underline{B^l_{\lambda}}\right)\underline{E^j}\prec D\underline{\bf N}. \left(\sum_{l\in \text{Supp}(\lambda)}\underline{B^l_{\lambda}}\right).
$$
As in \re{maj-Bb}, we have
$$
\left\|\sum_{\lambda}\underline{\frac{1}{b_{\lambda}}P_{\lambda}\left(\frac{\tilde B_{\lambda}}{b_{\lambda}}\right)}\right\|\leq \sup_{\lambda}\frac{1}{(|\lambda|-\underline{c_{\lambda}})^2}\left\|  D\underline{\bf N}. \left(\sum_{\lambda}\sum_{l\in \text{Supp}(\lambda)}\underline{B^l_{\lambda}}\right)\right\|.
$$
Hence, for $\|z\|\leq r$, 
\beq\label{estim2}
\left\|\underline{\sum_{\lambda}\frac{1}{b_{\lambda}}P_{\lambda}\left(\frac{\tilde B_{\lambda}}{b_{\lambda}}\right)}\right\|\leq 2\epsilon.
\eeq
Collecting estimates \re{estim1} and \re{estim2}, we obtain
$$
\sup_{\norma z\leq r} \norma{\sum_{\lambda}\underline{\frac{\tilde
  B_{\lambda}}{b_{\lambda}}-\frac{1}{b_{\lambda}}P_{\lambda}\left(\frac{\tilde
    B_{\lambda}}{b_{\lambda}}\right)}  }\leq 4\epsilon\ ,
$$
and remarking that, for functions of class $\cN$ the projector
$J^{2m}$ does not increase the norm, one gets
\beq\label{estimU}
\sup_{\|z\|\leq r}\|\underline{U}(z)\|\leq 4\epsilon.
\eeq
\end{proof}

\subsection{Flow of normally analytic vector fields}\label{s.flow}
In this section we study the flow $\Phi^t$ of a vector field
$U\in\cN_{r}$. In particular we will prove the following Lemma

\begin{lemma}
  \label{flow}
Assume that $U\in\cN_r$ for some $r>0$ fulfills $\epsilon:=\norma{U}_r<\frac{\delta}{4e}$ and let $\bF\in\cNF_r $ and $\delta<r$. 
Then the family
$\phis\bF\equiv\{\phis F^i\}_{i\geq 1}$ is summable
normally analytic and, defining $S^i:=\phis F^i-F^i$ and $\tilde
S^i:=\phis F^i-F^i-[U,F^i]$, one has
\begin{equation}
  \label{resto}
\norma{\underline{\bS}}_{r-\delta}\leq
\frac{4}{\delta}\norma{\bF}_r\epsilon\ ,\quad
\norma{\underline{\tilde \bS}}_{r-\delta}\leq
\frac{8e}{\delta^2}\norma{\bF}_r\epsilon^2
\end{equation}
\end{lemma}
\proof To start with, we remark that, since $\sup_{\norma z< r}\normap
       {U(z)}\leq \norma{\underline U}_r$, $\forall \left|t\right|\leq
       1$, one has
       \begin{equation*}
        \norma{\Phi^t(z)-z}=\norma{\int_0^tU(\Phi^s(z))ds}\leq
        \norma{\int_0^tU(\Phi^s(z))ds}_+ \leq \delta 
       \end{equation*}
and therefore $z\in B_{r-\delta}$ implies $\Phi^t(z)\in B_r$ i.e.
$\Phi^t(B_{r-\delta})\subset B_r$ ($B_r$ denoting
       the ball in $H$ of radius $r$ centered at zero). Thus the flow
       is well defined and analytic at least up to $|t|=1$. By Taylor
       expanding in $t$ at $t=0$, one has
       \begin{equation}
         \label{taylor}
(\Phi^{-t})^*F^i=\sum_{k\geq0}\frac{t^kAd_U^k}{k!}F^i\ ,
       \end{equation}

where $Ad_UG:=[U,G]$. To estimate this family remark first that
$$
\underline{Ad_{U}F^i}\preceq \underline{DU\, F^i}+\underline{DF^i\, U}=:\aad
\underline F^i\ .
$$
Summing over $i$ one gets
$$
\sum_{i}\underline{Ad_UF^i}\preceq \aad\underline\bF\ ,
$$
and, by induction on $k$
$$
\sum_{i}\underline{Ad_U^kF^i}\preceq \aad^k\underline\bF\ .
$$
Thus we have
\begin{equation}
  \label{tg}
\sum_{i}\left(\underline{\phis F^i}-\underline F^i\right)\preceq
\sum_{k\geq 1}\frac1{k!} \aad^k\underline\bF\ .
\end{equation}
In order to estimate the r.h.s. remark first that, for any family
$\bG\in\cNF_{r-\delta-\delta_1}$ (for some $\delta,\delta_1\geq 0$),
we have, by Cauchy estimate
\begin{equation}
  \label{aa}
\norma{\aad\underline\bG}_{r-\delta-\delta_1-\delta_2}\leq
\frac{2}{\delta_2}\norma{\underline
  U}_r\norma{\underline{\bG}}_{r-\delta-\delta_1} \ .
\end{equation}
Fix now some $k\geq0$, define $\delta':=\delta/k$ and look for
constants $C^{(k)}_l$, $0\leq l\leq k$ s.t.
$$
\norma{\aad^l\underline\bF}_{r-l\delta'}\leq C^{(k)}_l \ .
$$
Of course, by \eqref{aa} they can be recursively defined by
$$
C^{(k)}_l={\frac{2}{\delta'}}C^{(k)}_{l-1}\norma{\underline U}_r\ ,\quad
C^{(k)}_0:= \norma{\underline\bF}_r\ ,
$$
which gives
$$
C^{(k)}_l=\left( \frac{2}{\delta'}\norma{\underline U}_r \right)^l
\norma{\underline \bF}_r\ ;
$$
taking $l=k$ this produces an estimate of the general term of the
r.h.s. of \eqref{tg}:
\begin{equation}
\label{aa.1}
\norma{\frac{\aad^k\underline\bF}{k!}}_{r-\delta}\leq
\frac{k^k}{k!}\norma{\underline\bF}_r \left( \frac{2}{\delta}\norma{\underline
  U}_r \right)^k \leq \frac{\norma{\underline\bF}_r  }{e} \left( \frac{2e}{\delta}\norma{\underline
  U}_r \right)^k\ ,
\end{equation}
where we used $k!\geq k^k e^{-k+1}$. Summing over $k\geq 1$ or
$k\geq 2$, one gets the thesis. \qed

Although the family $\bE$ is not summably normally analytic, its
composition with the flow has the following remarkable property.
\begin{lemma}
  \label{flow.E}
Assume that $U\in\cN_r$ for some $r>0$ fulfils 
$\epsilon:=\norma{U}_r<\frac{\delta}{8e}$ with $0<\delta<r$; then the family
$\bT\equiv\{\phis E^i-E^i-[U,E^i]\}_{i\geq 1}$ is summably normally analytic
and one has
\begin{equation}
  \label{resto.E}
\norma{\underline{\bT}}_{r-\delta}\leq
\frac{8r}{e\delta}  \left(\frac{4e\epsilon}{\delta}\right)\epsilon\ .
\end{equation}
\end{lemma}
\proof We proceed as in the proof of the previous Lemma except that we
compute explicitly the first term of the expansion \eqref{taylor}.

One has $D\underline{U}\,\underline{\bE}=(D\underline U(z))z$ and
$\sum_i(D\underline{E^i})\underline{U}=\underline U$, so we get (for
any $\delta'<r$),
\begin{equation}
  \label{est.E}
\norma{\underline{[U,\bE]}}_{r-\delta'}\leq\left(\frac{r}{\delta'}+1\right)\norma
{\underline U}_r\leq \frac{2r}{\delta'}\norma
{\underline U}_r\ .
\end{equation}
So, by \eqref{aa.1}, 
$$
\frac{1}{(k-1)!}\norma{\aad
  ^{k}\underline\bE}_{r-2\delta'}=\frac{1}{(k-1)!}\norma{\aad
  ^{k-1}\underline{[U,\bE]}}_{r-2\delta'}\leq
\frac{2r}{e\delta'}\left(\frac{2e}{\delta'}\epsilon\right)^{k-1}\epsilon \ ,
$$
thus
\begin{align*}
\sum_{k\geq 2} \frac{1}{k!}\norma{\aad
  ^{k}\underline\bE}_{r-2\delta'}\leq \sum_{k\geq 2} \frac{1}{k}
\frac{2r}{e\delta'}\left(\frac{2e}{\delta'}\epsilon\right)^{k-1}\epsilon\leq
\frac{2r}{e\delta'}\epsilon \sum_{k\geq
  1}\left(\frac{2e}{\delta'}\epsilon\right)^{k}\leq
\frac{2r}{e\delta'}\epsilon 2 \left(\frac{2e}{\delta'}\epsilon\right)\ .
\end{align*}
Taking $\delta'=\delta/2$ one gets the thesis.\qed

\subsection{Iteration}
We use $U$ to generate a change of variables which is the time 1 flow,
$\Phi$ of the system $\dot z=U(z)$. We have 
\begin{eqnarray}
\Phi_*^{-1}X_i &= & X_i+[-U,X_i]+O(2m+1)\nonumber\\
&=& NF^i_{\leq m}+B^i_{nres}+B^i_0+[NF^i_{\leq m}, -U]+O(2m+1)\nonumber\\
&=& \underbrace{NF^i_{\leq m}+B^i_0}_{=:NF^i_{\leq 2m}}+O(2m+1).\label{equ-conj0}
\end{eqnarray}
By assumption, $B_0^i=\sum_{j\geq 1}\tilde a_{i,j}E^j$ where $\tilde a_{i,j}$ are polynomials of degree $\leq 2m-1$ that are common first integrals of $\bf E$.

By assumption, we have ${\bf X}={\bf E}+ {\bf F}$ and there exists
$c_0$ such that 
\begin{equation}\label{cond-const}
\norma{\underline{\bF}}_{r_0}\leq c_0r_0^2
\end{equation} for some small parameter $r_0$. We also fix two large constants $c_1$
and $b{\geq 1}$ (we will track the dependence of everything on such
constants). Their precise value will be decided along the procedure.

We
denote $m:=2^k$, $k\geq 0$ then the sequences we are interested in are defined
by
\begin{align}
\label{seque}
q_m&:=m^{-\frac{b}{m}}\ ,\quad m=2^k\ ,\quad k\geq 0
\\
\epsilon_k&:=\frac{\epsilon_0}{4^k}\ ,\quad k\geq0
\\
\epsilon_0&=c_0r_0^{2}
\\
\delta_0&:=\frac{r_0}{2}\ ,
\\
\delta&:=\frac{1}{c_1}r_0\ ,\quad  \delta_k:=\frac{\delta}{4^k}\ ,\quad
k\geq 1
\\
\quad r_1&:=\frac{1}{4}(r_0-\delta_0)=\frac{1}{8}r_0\ ,\quad
r_{k+1}:=q_{2^k}\left(r_k-\delta_k\right)\ ,\quad k\geq 1\ .
\end{align}
In the appendix we will prove that the following properties hold
\begin{align}
\label{dk}
{d_k:=}\prod_{l=0}^{k-1}q_{2^l}&=\frac{1}{4^{b\left(1-\frac{k+1}{2^k}\right)}}\geq
4^{-b}\ ,
\\
\label{rk}
r_k&\geq\frac{1}{4^b}r_1-\frac{\delta}{3}=\frac{1}{4^b}\left(\frac{1}{8}
-\frac{4^b}{3}\frac{1}{c_1}
\right)r_0  =:r_{\infty}\geq\frac{r_0}{4^{b+2}}\ ,
\end{align}
provided $c_1\geq 4^{b+2}/3$. Actually we take
\begin{equation}
\label{c.1}
c_1=\frac{4^{b+2}}{3}\ ,\ \Longrightarrow\ r_{\infty}=\frac{r_0}{4^{b+2}}\ .
\end{equation}
We will also prove that 
\begin{equation}
\label{eps}
\sum_{l=0}^{k-1}\epsilon_l\leq\frac{4}{3}\epsilon_0\ ,\quad
\sum_{l=0}^{k-1} \frac{\epsilon_l}{r_l-r_{l+1}}\leq
\frac{8}{7}\frac{\epsilon_0}{r_0
}+\frac{\epsilon_0}{r_{\infty}}2^{b/2}\ .
\end{equation}

Consider the following inequlities (with $m=2^k$) 
\begin{align}
\label{i.1}
\norma{\underline{\bR_{\geq m+1}}}_{r_k}&\leq \epsilon_k\ ,
\\
\label{i.2}
\norma{\underline{\bN_{\leq m}}}_{r_k}&\leq \left\{
  \begin{matrix}
    0&{\rm if} & k=0
    \\
    \sum_{l=0}^{k-1}\epsilon_l&{\rm if} & k\geq 1
  \end{matrix}
  \right.
\\
\label{i.3}
\sup_{\norma z\leq r_k}\left|D\underline{ \bN_{\leq m}}
(z)\right|_{\cB(H,H^+)}&\leq
\left\{
  \begin{matrix}
0&{\rm if} & k=0
\\
\sum_{l=0}^{k-1}\frac{\epsilon_l}{r_l-r_{l+1}}&{\rm if} & k\geq 1
  \end{matrix}
  \right.
\ .
\end{align}

\begin{lemma}\label{lemma-iter}
Assume
$$
b\geq 8+\frac{2\ln(3\cdot16)}{\ln 2}\ ,
$$
take 
\begin{align}
\label{sti.r}
c_1&=\frac{4^{b+2}}{3}\ ,
\\ r_0&<\min\left\{\sqrt{\frac{3}{8c_0}},{\frac{3}{136c_0};}\frac{1}{32 e  c_0c_1};\frac{1}{7\cdot  2^9\blur{e}c_0c_1^2};\frac{c_1{-\frac{1}{4}}}{c_0\left(2^{{4}}c_1{+2^9c_1^2}+4\frac{9}{7}+4\cdot
  2^{b/2}4^{b+2}\right)} \right\} \ .
\end{align}
Assume that the inequalities \eqref{i.1}-\eqref{i.3} hold with some
$k\geq 0$. Let $\Phi_m$ be the flow generated by $U_m$ defined by \re{Um}. It
conjugates the family $X^i_m= E^i+N^i_{\leq m}+R^i_{\geq m+1}$ to the
family $X^i_{2m}:=E^i+N^i_{\leq 2m}+R^i_{\geq 2m+1}$ and \eqref{i.1},
\eqref{i.2}, \eqref{i.3} hold for the new $N$ and $R$ with $k+1$ in
place of $k$.
\end{lemma}
\begin{proof}
First we define
\begin{equation}
  \label{n+}
N^i_{\leq 2m}:=N^i_{\leq m}+\left(J^{2m}R^i_{\geq m+1}\right)_{res}\ ,
\end{equation}
so that the estimate \eqref{i.2} immediately follows and the estimate
\eqref{i.3} follows from Cauchy inequality.

Then, an explicit computation gives
\begin{align}
  \label{r.1.1}
  R_{\geq 2m+1}^i&=\phim E^i-E^i-[U,E^i]
  \\
  \label{r.1.2}
  &+\phim N_{\leq m}^i-N_{\leq m}^i-[U,N_{\leq m}^i]
  \\
\label{r.2}
  &+\phim R_{\geq m+1}^i-R_{\geq m+1}^i
\\
\label{r.3}
&+\left(\id-J^{2m}\right)\left([U,E^i]+R^i_{\geq m+1}+[U,N^i_{\leq
    m}]\right)\ .
\end{align}
We remark that, as it can be seen by a qualitative analysis and we
will also see quantitatively, the largest contribution to the
estimate of the reminder term comes from the term $[U,E^i]$ in
\eqref{r.3}, followed (in size, but not in terms of order of
magnitude) by the term coming from $R_{\geq m+1}$ still in
\eqref{r.3}. All the other terms admit estimates which of higher
order.

 {Let us prove by induction on $k\geq 0$ estimates \re{i.1},\re{i.2} and \re{i.3}} 

\blur{\rv( We come to the precise estimate. To this end,we split the case $k=0$
(which corresponds to $m=1$) and the case $k\geq 1$.)}

For $k=0$, one has $N_{\leq 1}\equiv 0$  {and $\norma{\underline{\bf R}_{\geq 2}}_{r_0}\leq c_0r_0^2=\epsilon_0$. Hence, inequalities hold true for $k=0$. Assume that they hold for all $0\leq l\leq k$ and let us prove the inequality for $m=2^{k+1}$.}

 {Since $r_1$ and $\delta_0$ do not follow the induction definition of $r_k$ and $\delta_k$, we have to prove separatly the case $k=1$.}
Since $N_{\leq 1}\equiv 0$ then \eqref{r.1.2} is not
present, as well as the last term in \eqref{r.3}. Furthermore the
nonlinear cohomological equation reduces to the linear one, so $U$ can
be estimated using Remark \ref{rem-sd} with $c=1$ which gives
$$
\norma{U}_{r_0}\leq \epsilon_0\ .
$$ 
\blur{\rv(So, for $k=0$,} We have that, by \eqref{resto.E}, \eqref{resto}, the
families corresponding to \eqref{r.1.1} and \eqref{r.2} are estimated
by (with a little abuse of notation)
\begin{align*}
\norma{\underline{\eqref{r.1.1}}}_{r_0- \delta_0}&\leq
\frac{8r_0}{e\delta_0}\left(\frac{4e\epsilon_0}{\delta_0}\right)\epsilon_0
= 32\left(\frac{r_0}{\delta_0}\right)^2c_0r_0\epsilon_0 {=32\cdot 4 c_0r_0\epsilon_0}
\\
\norma{\underline{\eqref{r.2}}}_{r_0-\delta_0}&\leq
\frac{4}{\delta_0}\epsilon_0\epsilon_0  {=8 c_0r_0\epsilon_0}\ .
\end{align*}
Concerning \eqref{r.3}, by \eqref{est.E} we have
$$
\norma{\underline{[U,\bE]}}_{r_{ {0}}-\delta_{ {0}}}\leq \epsilon_0
\frac{2r_0}{\delta_0}= {4\epsilon_0}\ ,
$$
and thus
$$
\norma{\underline{\eqref{r.3}}}_{r_0- \delta_0}\leq\epsilon_0
\frac{2r_0}{\delta_0}+\epsilon_0\blur{<} {=5\epsilon_0}\ .
$$
It follows that
$$
\norma{\underline{\bR_{\geq 3}}}_{r_0-
  \delta_0}\leq  {((32\cdot 4 +8)c_0r_0 + 5)\epsilon_0 \leq 8\epsilon_0}
\blur{\left(3\frac{r_0}{\delta_0} + \frac{4}{\delta_0}
\epsilon_0+32\left(\frac{r_0}{\delta_0}\right)^2c_0r_0\right)\epsilon_0
\leq 4 \frac{r_0}{\delta_0}\epsilon_0} \ , 
$$
provided
\begin{equation}
  \label{r0}
r_0\leq {\frac{3}{136 c_0}}
\blur{\frac{1}{4c_0}\left(1+8\left(\frac{r_0}{\delta_0}\right)\right)^{-1}=\frac{17}{4c_0}} 
\ .
\end{equation}
From Lemma \ref{zero} it follows that 
$$
\norma{\underline{\bR_{\geq 3}}}_{\frac{1}{4}(r_0-
  \delta_0)}\leq \frac{1}{4^3} \norma{\underline{\bR_{\geq 3}}}_{r_0-
  \delta_0}\leq  \frac{1}{4^3} 4 \frac{r_0}{\delta_0}\epsilon_0
=\frac{1}{8}\epsilon_0<\frac{\epsilon_0}{4}=\epsilon_1\ .
$$
We also remark that, by Cauchy estimate, we have
$$
\norma{\underline{D\bN_{\leq 2}}}_{r_1}\leq \frac{1}{r_0-r_1}
\norma{\underline{\bN_{\leq 2}}}_{r_0}\leq \frac{1}{r_0-r_1}
\norma{\underline{\bR_{\geq 2}}}_{r_0} \leq
\frac{\epsilon_0}{r_0-r_1}\ .
$$
 {This concludes the proof of for case $k=1$.}

\vskip5pt 

Assume now $k\geq 1$. According to \rl{estim-nonlin-cohomo}, we have $\|\underline{U}(z)\|_{r_k}\leq 4\epsilon_k$ \blur{, so, we also have} {as soon as}
\begin{align}
\label{r0.1}
\frac{16e\epsilon_k}{\delta_k}=\frac{16e\epsilon_0}{\delta}=16ec_0c_1r_0
<&\frac{1}{2}\ \iff\ r_0<\frac{1}{32ec_0c_1}\\
\label{r0.2}
 {\frac{4}{3}\epsilon_0\leq} & \frac{1}{2}
\end{align} 
Hence, by \eqref{resto.E},
\eqref{resto}, the above families are estimated by
\begin{align*}
\norma{\underline{\eqref{r.1.1}}}_{r_k- \delta_k}&\leq
\frac{8r_k}{e\delta_k}\left(\frac{4e4\epsilon_k}{\delta_k}\right)4\epsilon_k
=  {2^9\left(\frac{\epsilon_0}{\delta}\right)^2r_k\leq}\blur{\rv(2^9\left(\frac{r_0}{\delta}\right)^2r_0c_0\epsilon_0=)} 2^9c_1^2r_0c_0\epsilon_0
\\
\norma{\underline{\eqref{r.1.2}}}_{r_k-\delta_k}&\leq
\frac{8e}{\delta_k^2}\frac{4}{3
}\epsilon_0(4\epsilon_k)^2=2^9c_0^2e\left(\frac{r_0}{\delta}\right)^2
r_0^2\epsilon_0 =\frac{2^9e}{ {3}}c_0^2c_1^2
r_0^2\epsilon_0 <2^9 c_0^2c_1^2
r_0^2\epsilon_0
\ ,
\\
\norma{\underline{\eqref{r.2}}}_{r_k-\delta_k}&\leq
\frac{4}{\delta_k}\epsilon_k4\epsilon_k=\frac{2^4c_0}{4^k}
\left(\frac{r_0}{\delta}\right)r_0
\epsilon_0 \leq {2^4c_0}r_0
\left(\frac{r_0}{\delta}\right)
\epsilon_0= {2^4c_0}r_0
c_1
\epsilon_0  \ .
\end{align*}
Concerning \eqref{r.3}, \blur{\rv(again)} by \eqref{est.E} we have
$$
\norma{\underline{[U,\bE]}}_{r_k-\delta_k}\leq 4\epsilon_k
\frac{2r_k}{\delta_k}\leq 8\left(\frac{r_0}{\delta}\right)\epsilon_0 =8c_1\epsilon_0 
$$
and, by \eqref{i.3}, \eqref{i.2}, \eqref{i.1} and
\eqref{eps}, we have
\begin{align*}
\norma{\underline{[U,\bN_{\leq m}]}}_{r_k-\delta_k}\leq
\norma{\underline{DU \bN_{\leq
      m}}}_{r_k-\delta_k}+\norma{\underline{D\bN_{\leq
      m}U}}_{r_k-\delta_k} \\ \leq
\frac{4\epsilon_k}{\delta_k}\frac{4}{3}\epsilon_0+4\epsilon_k\sum_{l=0}^{k-1}\frac{\epsilon_l}{r_l-r_{l+1}}\leq
\frac{\epsilon_0}{\delta}\frac{16}{3}\epsilon_0+
4\epsilon_k\left(\frac{8}{7}+2^{b/2}\frac{r_0}{r_{\infty}}\right)
\frac{\epsilon_0}{r_0}  
\\ =\left[\frac{16}{3}c_1+4\frac{8}{7}+4\cdot 2^{b/2}\frac{r_0}{r_\infty}
  \right] c_0r_0\epsilon_0\ . 
\end{align*}
Summing up we have
$$ \norma{\underline{[U,\bE]+[U,\bN_{\leq m}]+\bR_{\geq
      m+1}}}_{r_k-\delta_k}\leq\epsilon_0\left[8c_1 {+\frac{1}{4^k}}+c_0r_0
  \left(\frac{16}{3}c_1+4\frac{8}{7}+4\cdot
  2^{b/2}\frac{r_0}{r_\infty}\right) \right] \ ,
$$
and therefore the same estimate holds  for
$\norma{\underline{\eqref{r.3}}}_{r_k-\delta_k}$.
Summing up the different contributions, we have 
\begin{align}
  \label{rm}
\norma{\underline{\bR_{\geq 2m+1}}
}_{r_k-\delta_k}\leq
\epsilon_0\left[8c_1 {+\frac{1}{4^k}}+c_0r_0\left(\blur{2^{10}} {2^4}c_1 {+2^9c_1^2}+4\frac{8}{7}+4\cdot 2^{b/2}
  \frac{r_0}{r_\infty} +2^9c_0\blur{e}c_1^2r_0
\right)
  \right]
\ ,
\end{align}
which, provided 
\begin{equation}
\label{r0sti}
r_0<\frac{ {4}}{7\cdot 2^9c_0\blur{e}c_1^2}\ , \quad r_0<(c_1 {-\frac{1}{4})}
\left[c_0\left(\blur{2^{10}} {2^4}c_1 {+2^9c_1^2}+4\frac{9}{7}+4\cdot 2^{b/2}\frac{r_0}{r_\infty}\right)
  \right]^{-1}  
\ ,
\end{equation}
gives
\begin{equation}
  \label{nuovaR}
\norma{\underline{\bR_{\geq 2m+1}}
}_{r_k-\delta_k}\leq \blur{4\cdot} 9c_1\epsilon_0\ .
\end{equation}
ow, from Lemma \ref{zero}, since $b\geq 2$, one has
\begin{equation}
  \label{no-more-problem}
\norma{\underline{\bR_{\geq 2m+1}}
}_{r_{k+1}}\leq q_m^{2m+1}\blur{4\cdot} 9c_1\epsilon_0=3\cdot
4^{b+ {2}}2^{-bk\left(2+\frac{1}{2^k} \right)} {\epsilon_0}\ . 
\end{equation}
For $k=1$ (which corresponds to $m=2$), we have
$$
\norma{\underline{\bR_{\geq 5}}
}_{r_{2}}\leq 3\frac{4^{b+ {2}}}{2^{\frac{5}{2}b}}\epsilon_0
=\frac{3\cdot 4^ {2}}{2^{b/2}}\epsilon_0\leq \frac{\epsilon_0}{4^2}\ ,
$$
provided
\begin{align*}
\frac{3\cdot
  4^ {2}}{2^{b/2}}<\frac{1}{2^4}\ \iff\ \ln(3\cdot 4^{ {2}})<(\frac{b}{2}- 4
)\ln 2 \ ,
\end{align*}
which is equivalent to
\begin{equation}
\label{b.1}
b>8+\frac{2\ln(3\cdot 4^ {2})}{\ln2}\ .
\end{equation}
For $k\geq 2$ we have
$$
\norma{\underline{\bR_{\geq 2m+1}}
}_{r_{k+1}}\leq 3\cdot 4^{b+ {2}}2^{-2bk}\epsilon_0\leq
\frac{\epsilon_0}{4^{k+1}}\ ,
$$
provided 
$$
3\cdot 4^ {2}<4^{b(k-1)-(k+1)}\ \iff\ \frac{\ln(3\cdot 4^ {2}
  )}{\ln4}<b(k-1)-(k+1)\ \iff\ b>\frac{k+1}{k-1}+\frac{1}{k-1}\frac{\ln(3\cdot
  4^ {2})}{\ln4}\ ,
$$
which, since the r.h.s. is a decreasing function of $k$, is implied by
\begin{equation}
\label{b.2}
b>3+ \frac{\ln(3\cdot
  4^ {2})}{2\ln4}\ ,
\end{equation}
which in turn is implied by \eqref{b.1}.
\end{proof}

From Lemma \ref{lemma-iter}, by a completely standard argument, the
following Corollary follows

\begin{corollary}
The sequence of transformations $\{\Psi_k\}_{k\geq 1}$ defined by $\Psi_k:=\Phi_{2^{k-1}}^{-1}\circ \cdots\circ \Phi_{1}^{-1}$ converges to an analytic transformation $\Psi$ in a neighborhood of the origin and it conjugates the family $\{X^i\}_{i\geq 1}$ to a a family of normal forms $\{NF^i\}_{i\geq 1}$.
\end{corollary}
\appendix

\section{A technical Lemma}

\begin{lemma}
  \label{q}
Equation \eqref{dk} holds.
\end{lemma}
\proof Denote by $d_k$ the l.h.s. of \eqref{dk}, one has
\begin{align*}
d_k=\exp\left(\sum_{l=0}^{k-1}\ln
m^{-\frac{b}{m}}\right)=\exp\left(-\sum_{l=0}^{k-1}
\frac{b}{2^l}\ln 2^l\right)=\exp\left(-\frac{b\ln 2}{2}\sum_{l=0}^{k-1}
\frac{l}{2^{l-1}}\right)
\\
=\exp\left(-\frac{b\ln
  2}{2}4\left(1-\frac{k+1}{2^k}\right)\right) \ ,
\end{align*}
where we used the formula
$$
\sum_{l=0}^{k-1}
\frac{l}{2^{l-1}}=4\left(1-\frac{k+1}{2^k}\right)\ .
$$
Now, the result immediately follows.\qed

\begin{lemma}
  \label{rkappa}
Equation \eqref{rk} holds.
\end{lemma}
\proof We use the discrete analogue of the formula of the Duhamel
formula, namely we make the substitution $r_k=d_ks_k$, where $d_k$ was
defined in the proof of Lemma \ref{q}. One gets
\begin{align*}
r_{k+1}=d_{k+1}s_{k+1}=q_{2^k}d_ks_{k+1}=q_{2^k}(d_ks_k-\delta_k)
\end{align*}
and thus
$$
s_{k+1}=s_k-\frac{\delta_k}{d_k}\ ,\quad
s_1=\frac{r_1}{d_1}=r_1\ ,
$$
from which
$$
s_k=s_1-\sum_{l=1}^{k-1}\frac{\delta_l}{d_l}\ .
$$ Now, one has
$$
\sum_{l=1}^{k-1}\frac{\delta_l}{d_l}=\sum_{l=1}^{k-1}\frac{\delta}{4^l}
\frac{4^b}{4^{b\frac{l+1}{2^l}}}\leq \sum_{l=1}^{k-1}\frac{\delta}{4^l}
{4^b}=\frac{4^b}{3}\delta \ .
$$
Thus,
$$
r_k\geq d_k\left(\frac{r_1}{d_1}-  \frac{4^b}{3}\delta \right)\ .
$$

\qed

\begin{lemma}
  \label{epsilon}
Equation \eqref{eps} holds.
\end{lemma}
\proof
The first inequality is trivial. We discuss the second one. Using the
definition of $r_{k+1}$, we have
\begin{equation}
  \label{epsilon.1}
\frac{\epsilon_k}{r_k-r_{k+1}}=\frac{\epsilon_k}{r_k(1-q_{2^k})+
  q_{2^k}\delta_k}
\leq\frac{\epsilon_k}{r_k(1-q_{2^k})} \ ;
\end{equation}
now, one has
$$
1-q_m=1-\exp\left(-\frac{b}{m}\ln m\right)\ ,
$$
which is of the form $1-e^{-x}$ with $x$ varying from $0$ to
$\frac{b}{2}\ln2$. Remarking that in an interval $[0,x_0]$ one has 
$$
1-e^{-x}\geq e^{-x_0}x\ ,
$$
we get 
$$
1-q_m\geq 2^{-b/2}\left(\frac{b}{m}\ln m\right)=\frac{b}{2^{k+b/2}}\ln
2^k=\frac{k}{2^k} \frac{b}{2^{b/2}}\ln2\ ,
$$
and thus, for $k\geq 1$, 
\begin{align*}
\frac{\epsilon_k}{r_k(1-q_{2^k})}\leq \frac{\epsilon_0}{r_{\infty}}
\frac{2^{b/2}}{b\ln2}\frac{2^k}{k}\frac{1}{4^k}=\frac{\epsilon_0}{r_{\infty}}
\frac{2^{b/2}}{b\ln2}\frac{1}{k2^k}\ .
\end{align*}
Now one has
\begin{align*}
\sum_{k\geq1}\frac{x^k}{k}=\sum_{k\geq1}\int_0^xy^{k-1}dy=\sum_{k\geq0
}\int_0^xy^{k}dy=\int_0^x\frac{1}{1-y}dy=\left[-\ln\left|1-y\right|\right]_0^x
=-\ln\left|1-x\right|\ ,
\end{align*}
which, for $x=1/2$, gives
$$
\sum_{k\geq 1}\frac{1}{k2^k}=\ln2\ ,
$$
and thus
$$
\sum_{l\geq  {2}}\frac{\epsilon_l}{ {r_{l-1}-r_l}}\leq
\frac{\epsilon_0}{r_{\infty}} 2^{b/2}\ .
$$ adding the  {first}\blur{zero-th} term, namely
$\frac{8}{7}\frac{\epsilon_0}{r_0}$, one gets the thesis immediately
follows. \qed

\begin{lemma}
  \label{zero}
Let $\bF$ be a summable normally analytic vector fiels with $F^i$ having a
zero of order $m$ at the origin for all $i$. Let $0<\alpha\leq 1$,
then
\begin{equation}
  \label{alpha.1}
\norma{\underline{\bF}}_{\alpha r}\leq \alpha^m\norma{\underline{\bF}}_{r}\ .
\end{equation}
\end{lemma}
\proof Consider the function
$\underline\bF(z)=\sum_{Q,i}F_{Q,i}z^Q\ve_i$; since all the
coefficients are positive one has, for any $i$,
$$
\sum_{Q}F_{Q,i}(\alpha z)^Q=\alpha^m\sum_{Q,i}\alpha^{|Q|-m}F_{Q,i}z^Q\leq \alpha^m\sum_{Q}F_{Q,i}z^Q
$$
Thus one gets
$$
\sup_{\norma z\leq \alpha r}\normap{\underline \bF(z)}=\sup_{\norma
  z\leq r}\normap{\underline \bF(\alpha z)} \leq \alpha^m \sup_{\norma
  z\leq  r}\normap{\underline \bF(z)}\  .
$$
\qed

\bibliographystyle{alpha}
\bibliography{biblio}

\def\cprime{$'$}

\end{document}